\documentclass[letterpaper, 10 pt, conference]{ieeeconf}  % Comment this line out
\IEEEoverridecommandlockouts                              % This command is only
\usepackage{graphicx}
\usepackage{amssymb,amsmath,dsfont,amsfonts}
\usepackage{epstopdf}
\usepackage{tikz}
\usepackage{cite}
\usepackage{color}
\usepackage{boxedminipage}

\usepackage{subfigure}

\usepackage{algorithmic}
\usepackage{algorithm}

\newtheorem{theorem}{Theorem}
\newtheorem{lemma}{Lemma}

\newtheorem{remark}{Remark}

\newcounter{l1}
\newcommand{\barablist}{\begin{list}{\arabic{l1}}{\usecounter{l1}}}

\usepackage{acronym}
\acrodef{WLS}{Weighted Least Squares}

\title{An Abrupt Change Detection Heuristic\\ with Applications to Cyber Data Attacks on Power System$\text{s}^{\pi}$}
\author{Borhan M. Sanandaj$\text{i,}^{\star}$ Eilyan Bita$\text{r,}^{\diamond}$ Kameshwar Pooll$\text{a,}^{\star}$ and Tyrone L. Vincen$\text{t}^{\circ}$% <-this % stops a space
\thanks{\textsuperscript{$\star$}Borhan M. Sanandaji and Kameshwar Poolla are with Department of Electrical Engineering and Computer Sciences, University of California, Berkeley, CA 94720, USA, email:
\{sanandaji, poolla\}@berkeley.edu.}
\thanks{\textsuperscript{$\diamond$}Eilyan Bitar is with the School of Electrical and Computer Engineering, Cornell University, Ithaca, NY 14850, USA, email: eyb5@cornell.edu.}
\thanks{\textsuperscript{$\circ$}Tyrone L. Vincent is with the Department of Electrical Engineering and Computer Science, Colorado School of Mines, Golden, CO 80401, USA, email: tvincent@mines.edu.}
\thanks{\textsuperscript{$\pi$}This work was supported in part by EPRI and CERTS under sub-award 09-206; PSERC S-52; NSF under Grants EECS-1129061, CPS-1239178, and CNS-1239274; the Republic of Singapore National Research Foundation through a grant to the Berkeley Education Alliance for Research in Singapore for the SinBerBEST Program; Robert Bosch LLC through its Bosch Energy Research Network funding program.}
}

\def\real{ \mathbb{R} }

\newcommand{\vc}[1]{\boldsymbol{#1}}

\newcommand{\Prob}[1]{{\bf P}\left\{#1\right\}}

\newcommand{\argmin}{\mathop{\rm argmin}}

\acrodef{CoM}{Concentration of Measure}
\acrodef{i.i.d.}{independent and identically distributed}
\acrodef{LTI}{Linear Time-Invariant}
\acrodef{LTV}{Linear Time-Variant}
\acrodef{LPV}{Linear Parameter-Varying}
\acrodef{RIP}{Restricted Isometry Property}
\acrodef{SVD}{Singular Value Decomposition}
\acrodef{CS}{Compressive Sensing}
\acrodef{DSP}{Digital Signal Processing}
\acrodef{CSI}{Compressive System Identification}
\acrodef{CTI}{Compressive Topology Identification}
\acrodef{CBD}{Compressive Binary Detection}
\acrodef{OMP}{Orthogonal Matching Pursuit}
\acrodef{MP}{Matching Pursuit}
\acrodef{ERC}{Exact Recovery Condition}
\acrodef{BOMP}{Block Orthogonal Matching Pursuit}
\acrodef{COMP}{Clustered Orthogonal Matching Pursuit}
\acrodef{CoSaMP}{Compressive Sampling Matching Pursuit}
\acrodef{KKT}{Karush-Kuhn-Tucker}

\acrodef{FIR}{Finite Impulse Response}
\acrodef{DFT}{Discrete Fourier Transform}
\acrodef{DCT}{Discrete Cosine Transform}
\acrodef{JL}{Johnson-Lindenstrauss}
\acrodef{ROC}{Receiver Operating Curve}
\acrodef{NP}{Neyman-Pearson}
\acrodef{ARX}{Auto Regressive with eXternal input} 
\acrodef{MISO}{Multi-Input Single-Output}
\acrodef{SISO}{Single-Input Single-Output}
\acrodef{MIMO}{Multi-Input Multi-Output}
\acrodef{BP}{Basis Pursuit}
\acrodef{LASSO}{Least Absolute Shrinkage and Selection Operator}
\acrodef{GLASSO}{Group LASSO}
\acrodef{NNG}{Non-Negative Garrote}
\acrodef{LARS}{Least Angle Regression}
\acrodef{I/O}{Input/Output}
\acrodef{SCADA}{Supervisory Control and Data Acquisition}
\acrodef{PMU}{Phasor Measurement Unit}

\usepackage{color}

\newcommand{\cut}[1]{}

\begin{document}

\maketitle
\thispagestyle{empty}
\pagestyle{empty}

%%%%%%%%%%%%%%%%%%%%%%%%%%%%%%%%%%%%%%%%%%%%%%%%%%%%%%%%%%%%%%%%%%%%%%%%%%%%%%%%
\begin{abstract}
We present an analysis of a heuristic for abrupt change detection of systems with bounded state variations. The proposed analysis is based on the Singular Value Decomposition (SVD) of a history matrix built from system observations. We show that monitoring the largest singular value of the history matrix can be used as a heuristic for detecting abrupt changes in the system outputs. We provide sufficient detectability conditions for the proposed heuristic. As an application, we consider detecting malicious cyber data attacks on power systems and test our proposed heuristic on the IEEE 39-bus testbed.      
\end{abstract}

\acresetall
%%%%%%%%%%%%%%%%%%%%%%%%%%%%%%%%%%%%%%%%%%%%%%%%%%%%%%%%%%%%%%%%%%%%%%%%%%%%%%%%
%%%%%%%%%%%%%%%%%%%%%%%%%%%%%%%%%%%%%%%%%%%%%%%%%%%%%%%%%%%%%%%%%%%%%%%%%%%%%%%%
\section{Introduction}
Fault detection and supervisory control are essential to ensure that a dynamical system is operating in normal conditions. These monitoring mechanisms are of higher importance for critical systems such as power systems. 
Any propagation of faults in a power system may have severe consequences in the electricity generation, transmission, or distribution. To this end, \ac{SCADA} systems are designed for controlling and monitoring different parts of a power grid. Traditionally, within \ac{SCADA} or other conventional supervisory control and monitoring centers, the state of the system under study is estimated at every sample time. The condition of the system is then tested by monitoring a metric based on the estimated state. An abrupt change in that metric is an indicator of the occurrence of some malfunctioning in the system dynamics.

%%%%%%%%%%%%%%%%%%%%%%%%%%%%%%%%%%%%%%%%%%%%%%%%%%%%%%%%%%%%%%%%%%%%%%%%%
\subsection{Designated Data Attacks}
In power systems, as an example, changes of the system dynamics have been traditionally considered as a result of meter aging and malfunctioning, electrical breakdown, or natural causes such as storm, lightening, etc. However, such changes might be the result of a~\emph{designated} cyber data attack to the system. In particular, with the emergence of smart grids and its smart hardware and software components such as smart meters, \acp{PMU}, intelligent control devices, etc., power systems (and other similar large-scale dynamical systems) are more vulnerable to such malicious data attacks. 
In fact, it has been recently shown that an attacker can design attacks that do not appear in the detection metrics and can pass conventional detection algorithms.
Such attacks, namely called \emph{unobservable attacks}, require a careful compromise of meter readings by the attacker. Altogether, these have motivated a great amount of research to address cyber data attack detection within smart grids.
%
%%%%%%%%%%%%%%%%%%%%%%%%%%%%%%%%%%%%%%%%%%%%%%%%%%%%%%%%%%%%%%%%%%%%%%%%%%%%%%%%
%%%%%%%%%%%%%%%%%%%%%%%%%%%%%%%%%%%%%%%%%%%%%%%%%%%%%%%%%%%%%%%%%%%%%%%%%%%%%%%%
\subsection{Related Work}

Recently, Liu et al.~\cite{liu2011false} considered scenarios in which an attacker designs attacks carefully such that the conventional bad data detection algorithms are not capable of detecting them. Inspired by their work, many other papers targeted this problem~\cite{liu2011false,fawzi2011secure,pasqualetti2011cyber,giani2011smart,kosut2010malicious,kim2011strategic,gorinevsky2009estimation,bobba2010detecting,sandberg2010security}. An adversary attack has an impact on the real-time and day-ahead electricity markets. Such situations have been studied by~\cite{jia2011malicious,xie2010false}, among others.

Kosut et al.~\cite{kosut2010malicious} assume a Bayesian model on the state variables and consider a binary detection problem. In particular, they assume that the state variables have a zero-mean Gaussian distribution.  
Fawzi et al.~\cite{fawzi2011secure} impose a linear state-space representation on the power system state evolution and propose a decoder that corrects for the compromised meters. In their plant model, they assume they know the state transition and measurement matrices.

%%%%%%%%%%%%%%%%%%%%%%%%%%%%%%%%%%%%%%%%%%%%%%%%%%%%%%%%%%%%%%%%%%%%%%%%%%%%%%%%
\subsection{Main Contributions}
In this paper, we assume no a-priori distribution on the attack vector. We assume that the state variations (under normal conditions) are unknown but bounded within an $\ell_2$-norm. The time of attack (modeled as an abrupt change added to the unknown systems dynamics under normal conditions) and its magnitude is unknown to us as well. We present a heuristic for detecting such changes. The proposed heuristic is based on the \ac{SVD} of a history matrix built form system observations. We show that monitoring the  largest singular value of the history matrix is a good heuristic for detecting abrupt changes in the system outputs. In particular, we provide sufficient detectability conditions for the proposed heuristic. While the results of this paper can be applied to any system with a similar linear model with bounded state variations and generic faults, of our particular interest are power systems and unobservable attacks where such fault detection schemes play an important role in maintaining the safety and stability of the system.

%%%%%%%%%%%%%%%%%%%%%%%%%%%%%%%%%%%%%%%%%%%%%%%%%%%%%%%%%%%%%%%%%%%%%%%%%%%%%%%%
\subsection{Notation}
A column vector is shown as $\vc{x} \in \real^N$ with boldface letters. An element of a vector is shown as $x(i)$. A matrix is shown with a capital letter as $A \in \real^{M \times N}$. The elements of a matrix is shown as $A(i,j)$. The transpose and pseudo-inverse of $A$ are denoted by $A^T$ and $A^{\dagger}$, respectively. All variables are real-valued unless mentioned otherwise.

%%%%%%%%%%%%%%%%%%%%%%%%%%%%%%%%%%%%%%%%%%%%%%%%%%%%%%%%%%%%%%%%%%%%%%%%%%%%%%%%
%%%%%%%%%%%%%%%%%%%%%%%%%%%%%%%%%%%%%%%%%%%%%%%%%%%%%%%%%%%%%%%%%%%%%%%%%%%%%%%%
\section{Setup}
While the proposed analysis can be applied to any system with a linear model, we present our problem formulation based on the DC power flow model of a power system and a transmission network.
\subsection{Measurement Model}
Let $\vc{y}^t \in \mathbb{R}^{M}$ contain the injected power measurements of $n+1$ buses and line power measurements of $m$ branches of a transmission network at time $t$, where $M := m+n+1$.
Under a DC power flow model assumption over a finite time interval $t=t_i,\dots,t_f$, one can consider a linear relation between the measurements $\vc{y}^t$ and the power systems state vector $\vc{x}^t\in \real^N$ as: 

\begin{equation}
\vc{y}^t = H\vc{x}^t + \vc{e}^t, \quad t = t_i,\dots,t_f,
\label{eq:DCmodel-no-attack}
\end{equation}
where $\vc{x}^t \in \real^N$ is the state of the system at time $t$ containing the relative bus phase angles\footnote{The key state variables in a power grid contain bus voltage magnitudes and angles. However, in a DC power flow model the state variables are usually the bus voltage angles only.} and $\vc{e}^t$ contains the measurement noise. The matrix $H$ relates the state variables and the meter measurements and in general, is affected by the grid topology and link impedances. It is the job of the control center to construct $H$. In this paper, we assume $H$ is given and fixed over time and both the attacker and the control center have access to it. We assume the measurement noise is Gaussian with $\vc{e}^t \sim \mathcal{N}(0,\Lambda)$.
In order to incorporate the attack in the model, one can extend (\ref{eq:DCmodel-no-attack}) as:    
\begin{equation}
\vc{y}^t = H\vc{x}^t - \theta^t\vc{a} + \vc{e}^t, \quad t = t_i,\dots,t_f,
\label{eq:DCmodel-with-attack}
\end{equation}
where $\vc{a}$ is the attack vector and $\theta^t$ is an indicator variable. 

%%%%%%%%%%%%%%%%%%%%%%%%%%%%%%%%%%%%%%%%%%
%%%%%%%%%%%%%%%%%%%%%%%%%%%%%%%%%%%%%%%%%%
\subsection{Attack Model}
The attacker \emph{abruptly} changes the meter readings at time $t = t_a$ where $t_a$ is the time of attack. The indicator variable $\theta^t$ is defined as:
\begin{align} 
\theta^t = \left\{ \begin{array}{ll} 0, &  t < t_a, \\ 1, & t \geq t_a. \end{array} \right.
\label{eq:attack_model_1}
\end{align} 
%
%In a more delicate attack design, the attacker \emph{gradually} changes the meter readings. Formally, the indicator variable $\theta_t$ is defined as:
%%
%\begin{align} 
%\theta^t = \left\{ 
%\begin{array}{ll} 0, &  t < t_{a_i}, \\ 
%\frac{t-t_{a_i}+1}{t_{a_f}-t_{a_i}+1}, & t_{a_i} \leq t < t_{a_f},  \\
%1, & t \geq t_{a_f}, \end{array} 
%\right.
%\label{eq:attack_model_2}
%\end{align}
%% 
%where $t_{a_i}$ is the time an attack starts and $t_{a_f}$ is the time the the attack reaches it is final value, $\vc{a}$. Figure shows ...
%
\begin{remark}
In a more general setting, one can consider an arbitrary function $s(t,t_{a_i},t_{a_f})$ as the \emph{signature} of the attack where $t_{a_i}$ is when an attack starts and $t_{a_f}$ is when an attack reaches its final value. Apparently, there exists a trade-off for the attacker between the detectability of the attack (in this case, function $s$ should be smooth and gradually increasing rather than a step function) and the harmfulness of the attack (a step function has a larger harmful impact).
\hfill $\square$
\end{remark}

%%%%%%%%%%%%%%%%%%%%%%%%%%%%%%%%%%%%%%%%%%%%%%%%%%%%%%%%%%%%%%%%%%%%%%%%%%%%%%%%

\subsection{Systems with Bounded State Variations}

In this paper, we are interested in linear systems whose state variations are bounded within an $\ell_2$-norm ball. Formally, we consider systems with
\begin{equation}
\|\vc{x}^{t}-\vc{x}^{t_0}\|_2 \leq \gamma,
\label{eq:slowChanges}
\end{equation}
for any $t,t_0 \in \{t_i,\dots,t_f\}$ and for some $\gamma>0$. 
%%%%%%%%%%%%%%%%%%%%%%%%%%%%%%%%%%%%%%%%%%%%%%%%%%%%%%%%%%%%%%%%%%%%%%%%%%%%%%%%
\section{State Estimation}
\subsection{Before Attack}
Let's consider the attack model
(\ref{eq:DCmodel-with-attack})-(\ref{eq:attack_model_1}). Note that at any given time $t$ before the attack (i.e., $t_i \leq t < t_a$), we have $\vc{y}^t = H\vc{x}^t + \vc{e}^t$. 
A state estimate, $\widehat{\vc{x}}^t$, can be found by minimizing the cost function associated with a \ac{WLS} estimator as: 
\begin{equation*}
\widehat{\vc{x}}^t = \argmin_{\vc{x}^t} J(\vc{x}^t),
\end{equation*}
where
\[
J(\vc{x}^t):= (\vc{y}^t-H\vc{x}^t)^T\Lambda^{-1}(\vc{y}^t-H\vc{x}^t).
\]
It is trivial to find the minimizer of $J(\vc{x}^t)$. We have
\begin{equation}
\widehat{\vc{x}}^t = K\vc{y}^t,
\label{eq:state_WLS_estimate}
\end{equation}
where
\begin{equation}
K := (H^T\Lambda^{-1}H)^{-1}H^T\Lambda^{-1}.
\label{eq:pseudo-inverse-H}
\end{equation}
Substituting the measurement model $\vc{y}^t = H\vc{x}^t + \vc{e}^t$ in~(\ref{eq:state_WLS_estimate}), 
\[
\widehat{\vc{x}}^t= K\vc{y}^t = KH\vc{x}^t + K\vc{e}^t = \vc{x}^t + K\vc{e}^t,
\]
where we used the fact that $KH = I_N$.

%%%%%%%%%%%%%%%%%%%%%%%%%%%%%%%%%%%%%%%%%%%%%%%%%%%%%%%%%%%%%
%%%%%%%%%%%%%%%%%%%%%%%%%%%%%%%%%%%%%%%%%%%%%%%%%%%%%%%%%%%%%
\subsection{After the Attack}
After the attack, we have $\vc{y}^t_a := \vc{y}^t + \vc{a} = H\vc{x}^t + \vc{e}^t$. A \ac{WLS} estimate of the state under the attack can be found as: 

\begin{align}
\widehat{\vc{x}}_a^t &= K\vc{y}_a^t  = K\vc{y}^t + K\vc{a} \notag\\
&= \widehat{\vc{x}}^t+(H^T\Lambda^{-1}H)^{-1}H^T\Lambda^{-1}\vc{a}.
\label{eq:state-estimate-with-attack}
\end{align}
There has been a recent interest in the so-called \emph{unobservable} malicious data attacks on the power system~\cite{liu2011false}. From (\ref{eq:state-estimate-with-attack}), one can see that if
there exists a vector $\vc{c}\in \real^N$ such that
\[
\vc{a} = H\vc{c},
\]
then we have
\begin{equation}
\widehat{\vc{x}}_a^t = \widehat{\vc{x}}^t+\vc{c}
\label{eq:state_WLS_estimate_under_atatck}
\end{equation}
and consequently,
\[
\vc{y}_a^t-H\widehat{\vc{x}}_{\vc{a}}^t = \vc{y}^t-H\widehat{\vc{x}}^t.
\]
Thus, at each sample time $t$, if the residual $\vc{r}^t:=\vc{y}^t-H\widehat{\vc{x}}^t$ passes any detection metric (e.g., $\|\vc{y}^t-H\widehat{\vc{x}}^t\|_2$), the residual under attack $\vc{r}_a^t:=\vc{y}_a^t-H\widehat{\vc{x}}_a^t$ would pass that criteria, and consequently such an attack is unobservable from the view point of the control center. To this end, such detection algorithms are namely referred to as ``bad'' detection algorithms.

%%%%%%%%%%%%%%%%%%%%%%%%%%%%%%%%%%%%
%
%\subsection{State Estimate Model Under Unobservable Attacks}
%
%Combining the state estimates before and after attack (given in~(\ref{eq:state_WLS_estimate}) and~(\ref{eq:state_WLS_estimate_under_atatck})), we can consider a state estimate model over time as:
%%
%\begin{equation}
%\widehat{\vc{x}}^t = \vc{x}^t + K\vc{e}^t + \theta^t\vc{c},
%\label{eq:state_estiamtes}
%\end{equation}
%%
%where $\vc{x}^t$ is the true state at time $t$ and $\vc{c}$ is the change vector such that $\vc{a} = H\vc{c}$. 
%The attack vector $\vc{a}$ is assumed to be sparse. 
%%%%%%%%%%%%%%%%%%%%%%%%%%%%%%%%%%%%
%
%%%%%%%%%%%%%%%%%%%%%%%%%%%%%%%%%%%%
%
%%%%%%%%%%%%%%%%%%%%%%%%%%%%%%%%%%%%
%
\section{Analysis of an SVD-based Heuristic for Abrupt Change Detection}
% of Slowly Time-Varying Systems}

%In this section, we present an algorithm for detecting unobservable sparse attacks, i.e., attack vectors $\vc{a}\in \real^M$ with $\|\vc{a}\|_0 \leq S$ ($S\ll M$) and $\vc{a} = H\vc{c}$, for some $\vc{c} \in \real^N$. 
%
%Our goal is to leverage the idea of observable islands associated with unobservable sparse attacks. Note that all buses within an observable island have the same perceived voltage angle change, resulting in zero power flow perturbations on lines within any observable island. 
%
In this section, we analyze an \ac{SVD}-based heuristic which can be used for abrupt change detection of systems with bounded state variations. In our proposed approach, we collect a trace of the measurements over a finite time window.
%%%%%%%%%%%%%%%%%%%%%%%%%%%%%%%%%%%%
%
\subsection{History Matrix $\Delta^t$}
\label{subsec:historyMatrix}
Given the measurements $\vc{y}^t$ over a finite horizon of time, at any given time $t$ one can build a \emph{history matrix} that contains the changes of the measurements as:
\begin{equation}
\Delta^t =
\begin{bmatrix}
 ({\vc{y}}^t-{\vc{y}}^{t-1})^T \\
 ({\vc{y}}^t-{\vc{y}}^{t-2})^T \\
 \vdots \\
 ({\vc{y}}^t-{\vc{y}}^{t-w})^T
\end{bmatrix}^T \in \real^{M \times w},
\end{equation}
where $w$ is the size of the considered time window.
Define
\begin{equation*}
E^{t} :=
 \begin{bmatrix}
 {\vc{e}^{t}}^T \\
 {\vc{e}^{t}}^T \\
 \vdots \\
 {\vc{e}^{t}}^T
\end{bmatrix}^T \in \real^{M \times w}, \ \
G^{t} := 
\begin{bmatrix}
 {-\vc{e}^{{t}-1}}^T \\
 {-\vc{e}^{{t}-2}}^T \\
 \vdots \\
 {-\vc{e}^{{t}-w}}^T
\end{bmatrix}^T\in \real^{M \times w},
\end{equation*}
\begin{equation*}
X^{t} := 
\begin{bmatrix}
 {(\vc{x}^t-\vc{x}^{{t}-1}})^T \\
 {(\vc{x}^t-\vc{x}^{{t}-2}})^T \\
 \vdots \\
 {(\vc{x}^t-\vc{x}^{{t}-w}})^T
\end{bmatrix}^T\in \real^{N \times w}, \ \text{and} \ \ 
A := 
\begin{bmatrix}
 -\vc{a}^T \\
 -\vc{a}^T \\
 \vdots \\
 -\vc{a}^T
\end{bmatrix}^T,
\end{equation*}
for any $t$. It is trivial to see that $\Delta^t$ can be decomposed as
\begin{equation}
\Delta^t = E^t + G^t + HX^t, \ \ \ (\forall t < t_a).
\label{eq:delta_before_attack}
\end{equation}
At $t = t_a$ (i.e., when the attack happens), 
%
%\begin{align}
%\Delta^{t_a} &= 
%\begin{bmatrix}
% ({\vc{e}}^{t_a}-{\vc{e}}^{{t_a}-1})^T \\
% ({\vc{e}}^{t_a}-{\vc{e}}^{{t_a}-2})^T \\
% \vdots \\
% ({\vc{e}}^{t_a}-{\vc{e}}^{{t_a}-w})^T
%\end{bmatrix}^T
%+
%\begin{bmatrix}
% -\vc{a}^T \\
% -\vc{a}^T \\
% \vdots \\
% -\vc{a}^T
%\end{bmatrix}^T \notag \\
%&= E^{t_a} + G^{t_a} + A.
%\label{eq:delta_at_attack}
%\end{align}
\begin{equation}
\Delta^{t_a} = E^{t_a} + G^{t_a} + HX^{t_a} + A.
\label{eq:delta_at_attack}
\end{equation}
%
%\begin{equation}
%\Delta^{t_a} = E^{t_a} + G^{t_a} + A,
%\label{eq:history_matrix_decomposed}
%\end{equation}
%
%
%\begin{remark}
Similarly, at $t = t_a+1$, we have
\begin{equation}
\Delta^{t_a+1} = 
%\begin{bmatrix}
% ({\vc{e}}^{t_a+1}-{\vc{e}}^{t_a})^T \\
% ({\vc{e}}^{t_a+1}-{\vc{e}}^{t_a-1})^T \\
% \vdots \\
% ({\vc{e}}^{t_a+1}-{\vc{e}}^{t_a-(w-1)})^T
%\end{bmatrix}^T
E^{t_a+1} + G^{t_a+1} + HX^{t_a+1}
+
\begin{bmatrix}
 \vc{0}^T \\
 -\vc{a}^T \\
 \vdots \\
 -\vc{a}^T
\end{bmatrix}^T
\label{eq:delta_after_attack}
\end{equation}
and for $t \geq t_a+w$,
\begin{equation}
\Delta^t =
%\begin{bmatrix}
% ({\vc{e}}^t-{\vc{e}}^{t-1})^T \\
% ({\vc{e}}^t-{\vc{e}}^{t-2})^T \\
% \vdots \\
% ({\vc{e}}^t-{\vc{e}}^{t-w})^T
%\end{bmatrix}^T, 
E^{t} + G^{t} + HX^t,
\qquad (t\geq t_a+w).
\label{eq:delta_after_attack_window}
\end{equation} 
Note that the structure of the history matrix for $t\geq t_a+w$ in~(\ref{eq:delta_after_attack_window}) is similar to the one for $t<t_a$ given in~(\ref{eq:delta_before_attack}). 
%\end{remark}
%%%%%%%%%%%%%%%%%%%%%%%%%%%%%%%%%%%%%%%%%%%%%%%%%%%%%%%%%%%%%%%%%%%%%%%%%%
\subsection{Singular Value Analysis on $\Delta^t$}

Based on the structure of $\Delta^t$ and how it changes over time (before and after the attack), one can consider a heuristic for detecting abrupt changes in $\vc{y}^t$. While the rank of $\Delta^t$ (the number of non-zero singular values) does not change before and after the attack, the distribution of the singular values of $\Delta^t$ changes (due to the addition of a rank-1 matrix) after the attack. In particular, there exists a large jump in the largest singular value of the history matrix at the time of attack and afterwards.\footnote{While it is still noticeable, this jump starts to decrease at later times after the attack and vanishes at $t = t_a+w$.} In what follows, we monitor this jump and provide bounds on its magnitude before and after the attack.  
In particular, note that $E^{t_a}$ and $A$ are rank-1 matrices. Our goal is to exploit such a structure (rank-1 structure of $A$ and $E^{t_a}$) in evaluating the changes in the first singular value of $\Delta^t$. In order to keep the paper self-contained, we provide all required lemmas and theorems in proving the main theorems. 

Let $\sigma_i(\Delta^t)$ denote the $i$th singular value of $\Delta^t$. The following theorems present probability tail bounds on $\sigma_1(\Delta^t)$ (the largest singular value of $\Delta^t$). The first theorem shows that the largest singular value of $\Delta^t$ is bounded from above with exponentially high probability when $t<t_a$. 
\begin{theorem}
Let $\tau >0$ and $\epsilon >0$. Consider a linear system described by~(\ref{eq:DCmodel-with-attack}) and with bounded state variations as described by~(\ref{eq:slowChanges}). Let $M$ be the number of measurements and $w$ be the window size. Assume $\vc{a}$ be an unknown attack vector and $\vc{e^{t}}\sim\mathcal{N}(0,\nu^2)$. Let $\Delta^{t}$ and $G^{t}$ be defined as in section~\ref{subsec:historyMatrix}. Then, for $t<t_a$
\begin{equation*}
\Prob{\sigma_1(\Delta^{t}) \geq  \ell} \leq 2\exp{(-\frac{\tau^2}{2})} + \big((1+\epsilon)e^{-\epsilon}\big)^{M/2},
\end{equation*}
where 
\[
\ell := \nu\sqrt{w}\sqrt{M}(1+\epsilon)+\nu(\sqrt{M}+\sqrt{w}+\tau)+\gamma\sqrt{w}\|H\|.
\]
\label{theo:lowerBnd_sigma1_whenNoAttack}
\end{theorem}
\begin{proof}
See Appendix.
\end{proof}

%
%\begin{lemma}
%%
%Consider the measurement model~(\ref{eq:DCmodel-with-attack}) and the history matrix~$\Delta^{t_a}$ as decomposed in~(\ref{eq:history_matrix_decomposed}) with a given $w\geq 1$. Assume $\vc{a}$ be such that 
%$\sqrt{w}\|\vc{e}^{t_a}+\vc{a}\|_2 > \sigma_{1}(G^{t_a})$. Then
%%
%\begin{equation}
%\sigma_1(\Delta^{t_a}) \geq 
%\sqrt{w}\|\vc{e}^{t_a}+\vc{a}\|_2 - \sigma_{1}(G^{t_a}).
%\label{eq:delta_1_lower_bnd}
%\end{equation}
%\label{lem:bounds_on_svd_1}
%\end{lemma}
%%
%\begin{proof}
%See Appendix.
%\end{proof}
%
%%
%We are interested in bounds on the largest singular value of the history matrix, before and at the time attack. From (\ref{eq:delta_1_lower_bnd}) we have
%%
%\begin{equation}
%\sigma_1(\Delta^{t_a}) \geq \sqrt{w}\|\vc{e}^{t_a}+\vc{a}\|_2 - \sigma_{1}(G^{t_a}).
%\label{eq:delta_diff_lower}
%\end{equation}
%%
%On the other hand, from~(\ref{eq:delta_1_upper_bnd})
%\begin{align}
%&\sigma_1(\Delta^{t_a-1}) - \sigma_2(\Delta^{t_a-1}) 
%\leq \sigma_1(\Delta^{t_a-1}) \notag\\ 
%&\leq \sqrt{w}\|\vc{e}^{t_a-1}\|_2 +  \sigma_1(G^{t_a-1}).
%\label{eq:delta_diff_upper}
%\end{align}
%
%%
%Based on (\ref{eq:delta_diff_lower}) and (\ref{eq:delta_diff_upper}), one can consider a detectability criterion as defined in the following. 
%
The second theorem shows that $\sigma_1(\Delta^{t_a})$ is bounded from below with exponentially high probability.
\begin{theorem}
Let $\tau >0$ and $\epsilon >0$. Consider a linear system described by~(\ref{eq:DCmodel-with-attack}) and with bounded state variations as described by~(\ref{eq:slowChanges}). Let $M$ be the number of measurements and $w$ be the window length. Assume $\vc{a}$ be an unknown attack vector and $\vc{e^{t}}\sim\mathcal{N}(0,\nu^2)$. Assume $\|\vc{a}\|_2 \geq \|\vc{e}^{t_a}\|_2$. Let $\Delta^{t_a}$ and $G^{t_a}$ be defined as in section~\ref{subsec:historyMatrix}. Then,
\begin{equation*}
\Prob{\sigma_1(\Delta^{t_a}) \leq  u} \leq 2\exp{(-\frac{\tau^2}{2})} + \big((1+\epsilon)e^{-\epsilon}\big)^{M/2},
\end{equation*}
where 
\[
u=: \sqrt{w}\|\vc{a}\|_2 - \ell,
\]
\label{theo:upperBnd_sigma1_whenAttack}
\end{theorem}
and $\ell$ is as defined in Theorem~\ref{theo:lowerBnd_sigma1_whenNoAttack}.

\begin{proof}
See Appendix.
\end{proof}
\begin{remark}
Theorems~\ref{theo:lowerBnd_sigma1_whenNoAttack} and~\ref{theo:upperBnd_sigma1_whenAttack} provide probability tail bounds on $\sigma_1(\Delta^{t})$ before and at the time of attack, respectively. The results have a probabilistic notion with exponential bounds. When we say ``with high probability'' it refers to such exponential behavior. With reasonable choices of $\tau$ and $\epsilon$ for a given $M$, the probability term $2\exp{(-\frac{\tau^2}{2})} + \big((1+\epsilon)e^{-\epsilon}\big)^{M/2}$ can be pushed to be very close to $0$.
\hfill $\square$
\end{remark}

Based on these results, a detection rule can be considered as follows. For any given $\ell$ and $u$ such that $\ell < u$, if $\sigma_1(\Delta^t) < \ell$ for all $t_i\leq t<t_a$ and $\sigma_1(\Delta^{t_a}) > u$, then an attack has happened at time $t_a$. Based on this detection rule, we can derive the detection probability as follows.

\begin{theorem}
Let $\tau >0$ and $\epsilon >0$. Consider a linear system described by~(\ref{eq:DCmodel-with-attack}) and with bounded state variations as described by~(\ref{eq:slowChanges}). Let $M$ be the number of measurements and $w$ be the window length. Assume $\vc{a}$ be an unknown attack vector and $\vc{e^{t}}\sim\mathcal{N}(0,\nu^2)$. Let $\ell$ and $u$ be as defined in Theorems~\ref{theo:lowerBnd_sigma1_whenNoAttack} and~\ref{theo:upperBnd_sigma1_whenAttack}, respectively. Then, an attack $\vc{a}$ can be detected at $t_a$ with detection probability 
\begin{equation*}
\Prob{\text{detection}} \geq 1-2\big[2\exp{(-\frac{\tau^2}{2})} + \big((1+\epsilon)e^{-\epsilon}\big)^{M/2}\big].
\end{equation*}
if
\begin{equation*}
\|\vc{a}\|_2 > 2\big[\nu\sqrt{M}(1+\epsilon+
\frac{1}{\sqrt{w}}+\frac{1}{\sqrt{M}}+\frac{\tau}{\sqrt{M}\sqrt{w}})+\gamma\|H\|\big].
\end{equation*}
\label{theo:detectionProb}
\end{theorem}

\begin{proof}
See Appendix.
\end{proof}

\begin{remark}
Theorem~\ref{theo:detectionProb} provides a sufficient condition on $\|\vc{a}\|_2$ for detectability. The derived bound illustrates how different factors affect the detectability. For example, the larger the noise level (i.e., large $\nu$ value), the harder the detection. The number of measurements $M$ and the window size $w$ are also affecting the detectability. Detection of abrupt changes in systems with smaller state variations (i.e., smaller $\gamma$) is also easier as can be interpreted from this result.
\hfill $\square$
\end{remark}

\section{Case Study - IEEE 39-Bus Testbed}

\begin{figure}[tb]
\begin{minipage}[b]{\linewidth}
\begin{center}
(a)
 \includegraphics[width =\columnwidth]{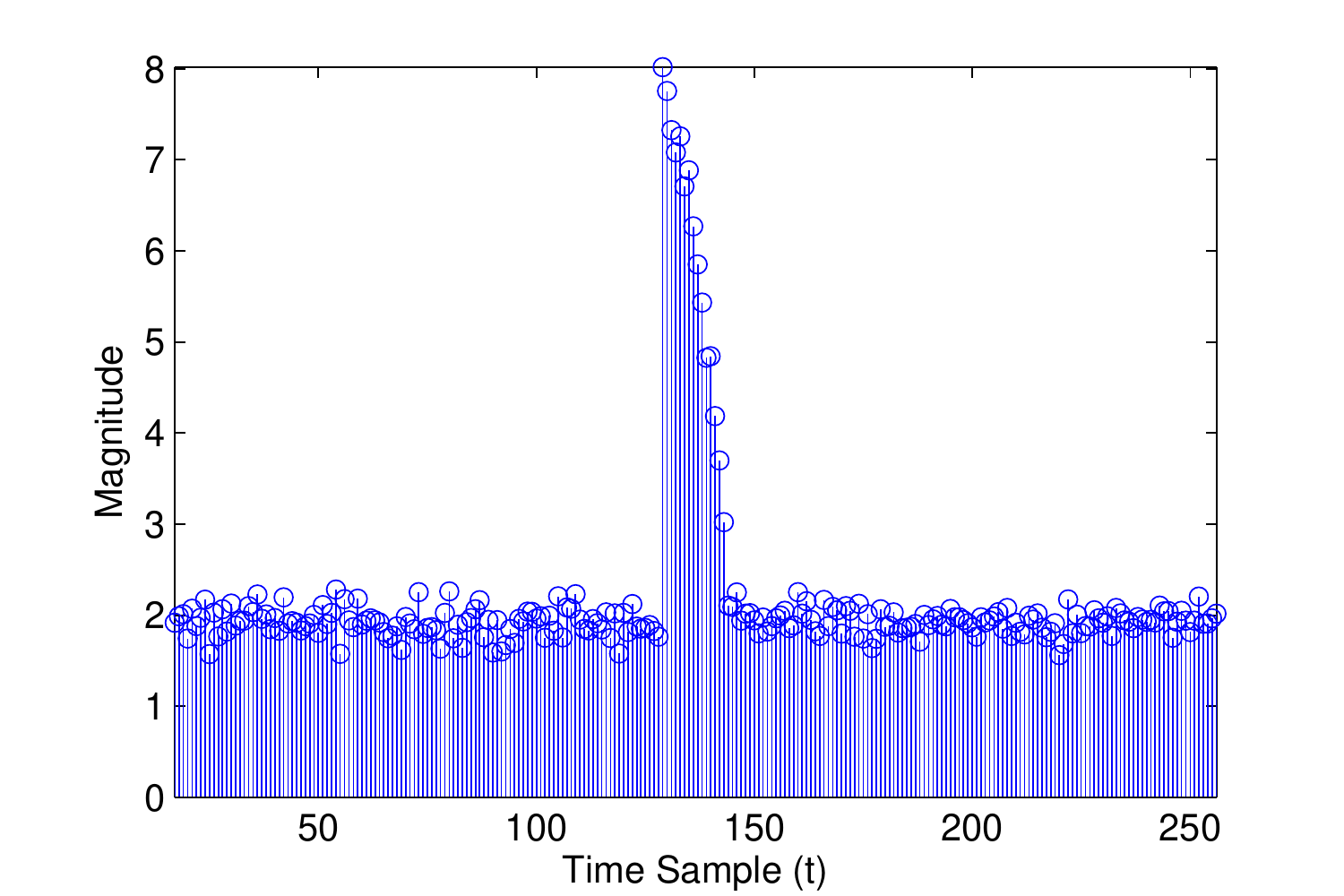}
 \end{center}
\end{minipage}
\begin{minipage}[b]{\linewidth}
\begin{center}
(b)
 \includegraphics[width = \columnwidth]{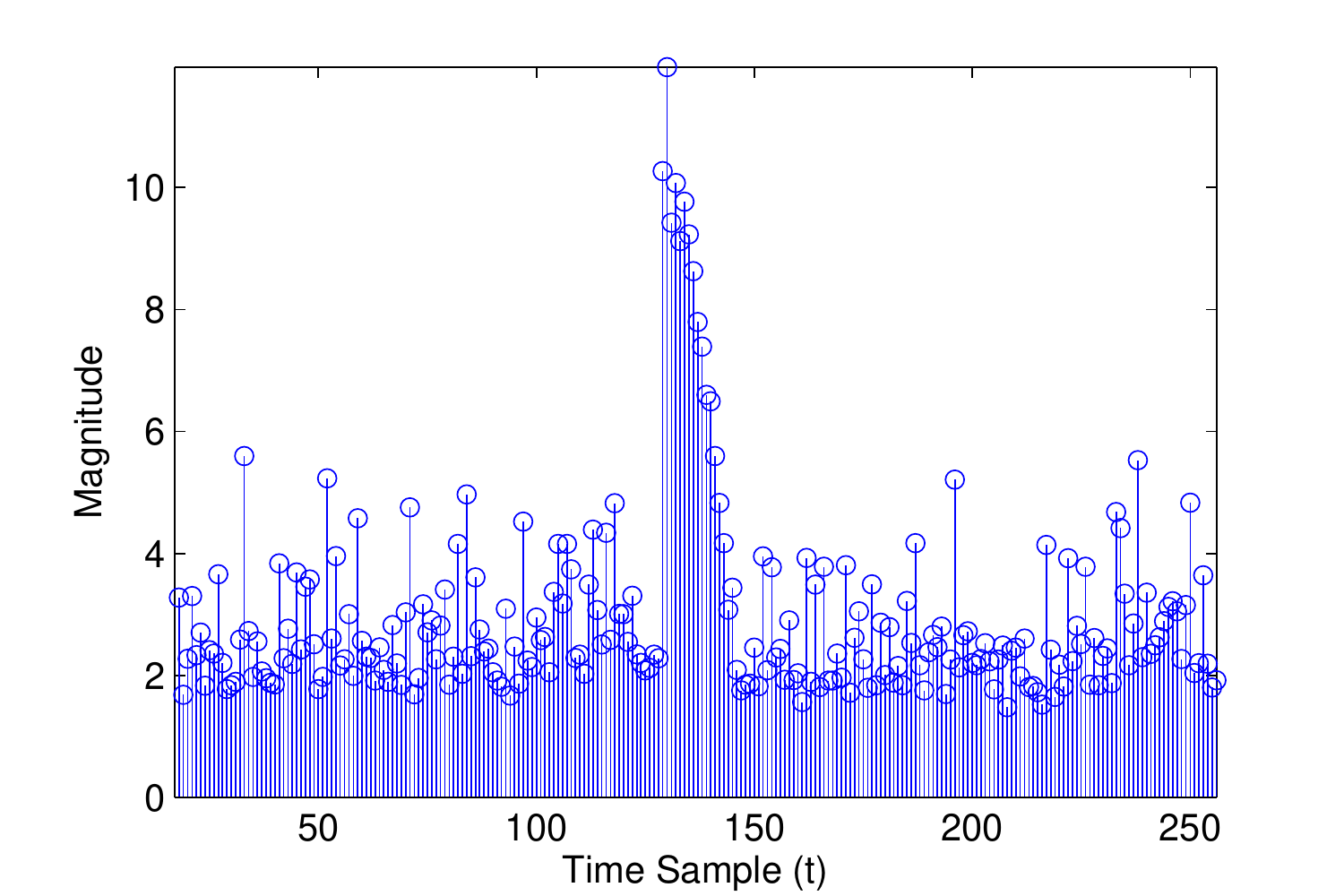}
  \end{center}
\end{minipage}
\caption{Illustration of how $\sigma_1(\Delta^t)$ changes over time when an attack happens.  In this example, an unobservable attack happens at $t=129$. A window size of $w=16$ and a noise level with $\nu = 0.05$ are considered. The history matrix is built based on (a) measurements and (b) state estimates.}
\label{fig:svd_test_plot1}
\end{figure}

In this section, we examine our proposed heuristic and detection condition on the IEEE 39-bus testbed~\cite{zimmerman2011matpower}. Consider a $4$-sparse unobservable attack happening at $t_a = 129$. Fig.~\ref{fig:svd_test_plot1} illustrates how $\sigma_1(\Delta^{t})$ evolves over time, where $M=85$ and $w=16$. 
As an alternative to measurements, one could build the history matrix based on the state estimates as given by (\ref{eq:state_WLS_estimate}). To this end, we consider two cases. Once we construct the history matrix based on the measurmenets (Fig.~\ref{fig:svd_test_plot1}(a)) and once based on the state estimates (Fig.\ref{fig:svd_test_plot1}(b)). As can be seen, the jump in $\sigma_1(\Delta^{t_a})$ is more distinguishable when the history matrix is built based on the measurements. In all of the simulations of this section, we assume $\gamma = 0$. However, similar results can be achieved for the case with $\gamma \neq 0$.

%Th detection probability is calculated as the percentage of the number of times an attack has been correctly detected over several noise levels (for each noise level, the detection rate is calculated over $1000$ trials). As expected, the detection performance decreases as the the noise level increases.
%%
Also, note that how $\sigma_1(\Delta^{t})$ starts to decrease at times after the attack (i.e., $t \geq 129$). In fact, this can be understood by looking at the structure of $\Delta^t$ at times after the attack. As shown in (\ref{eq:delta_after_attack}), the first column of the change matrix is zero at $t=t_a+1$. This makes the Frobenius norm (and consequently $\ell_2$-norm) of the history matrix smaller. This decrease in $\sigma_1(\Delta^{t})$ continues on until $t = t_a+w$. The effect of the attack disappears for $t \geq t_a+w$ when the characteristics of $\Delta^t$ are similar to the ones at $t<t_a$. 
\cut{
Fig.~\ref{fig:detection_performance} shows how $w$ and $q$ affect the detectability performance. As can be seen, the smaller $q$ and the larger $w$, the stronger the detection performance. One should note that our bound in Theorem~\ref{theo:bound_on_a} reflects these dependencies on $q$ and $w$. While increasing the window size would consequently increase the computation cost of our algorithm (mainly due to the \ac{SVD} step of the algorithm), one can overcome these issues by considering some recent methods for \ac{SVD} of large matrices~\cite{martinsson2011randomized} that significantly reduce the 
computational burden behind \ac{SVD}.

We would like to mention that the sufficient condition on $\|\vc{a}\|_2$ derived in Theorem~\ref{theo:bound_on_a} is reasonably tight. For instance, consider the detection result given in Fig.~\ref{fig:detection_performance}(a) where $q=2$, $w=64$, and $M=85$. For reasonable choices of $\epsilon=\sqrt{2}/2$ and $\tau=4$, Theorem~\ref{theo:bound_on_a} states that if $\|\vc{a}\|_2 \geq Q=58\nu$, with probability as least 0.9999 the attack vector is $2$-detectable. On the other hand, in this example $\|\vc{a}\|_2=2$ and based on the simulation results we have perfect detection performance for $\nu\leq 0.1$. Note that we have provided a sufficient condition for detectability and the gap between $\|\vc{a}\|_2=2$ and
$Q=5.8$ (when $\nu = 0.1$) is reasonable. 
}
%%%%%%%%%%%%%%%%%%%%%%%%%%%%%%%%%%%%%%%%%%%%%%%%%%%%%%

\begin{figure}[tb]
\begin{minipage}[b]{\linewidth}
\begin{center}
(a)
 \includegraphics[width =\columnwidth]{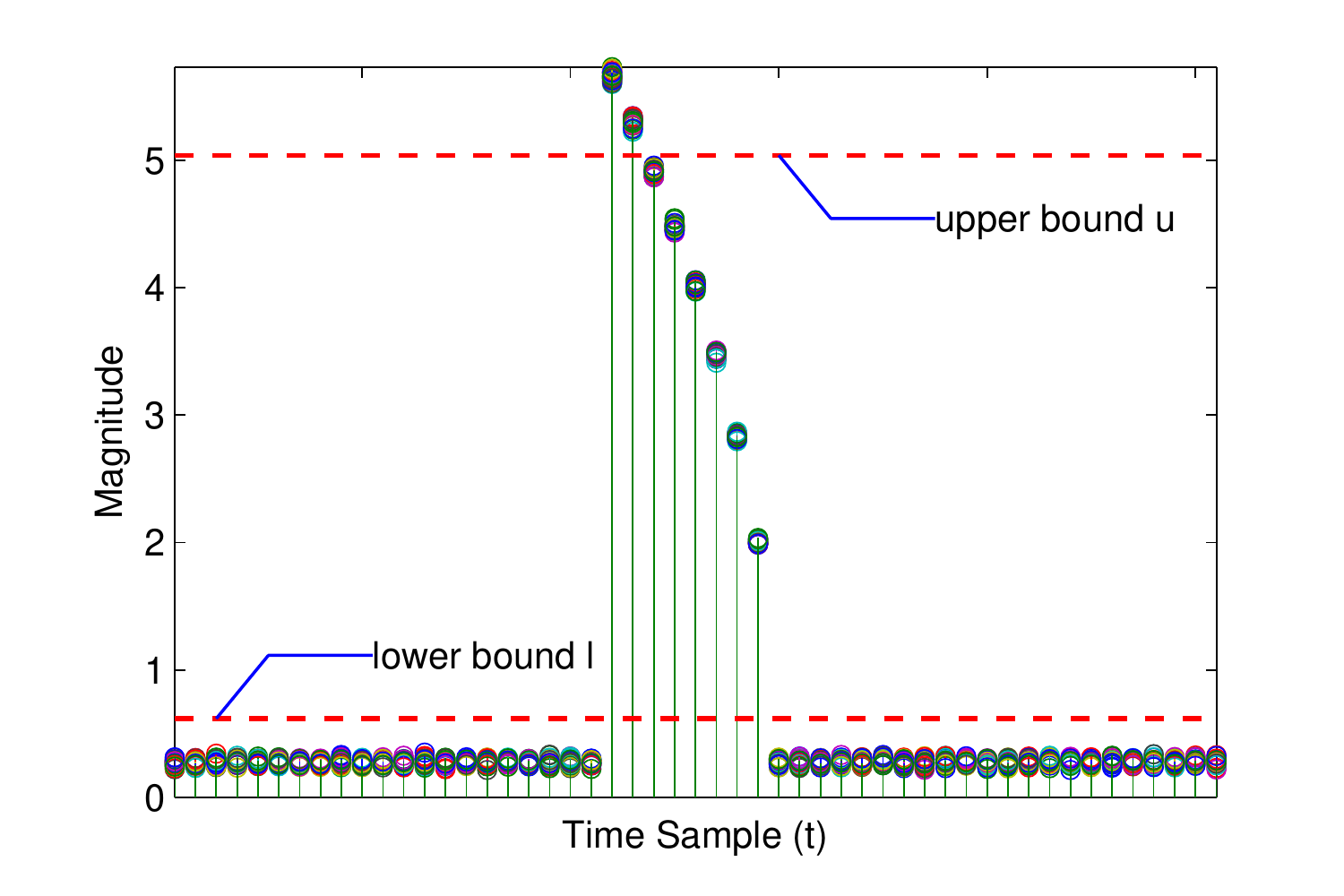}
 \end{center}
\end{minipage}
\begin{minipage}[b]{\linewidth}
\begin{center}
(b)
 \includegraphics[width = \columnwidth]{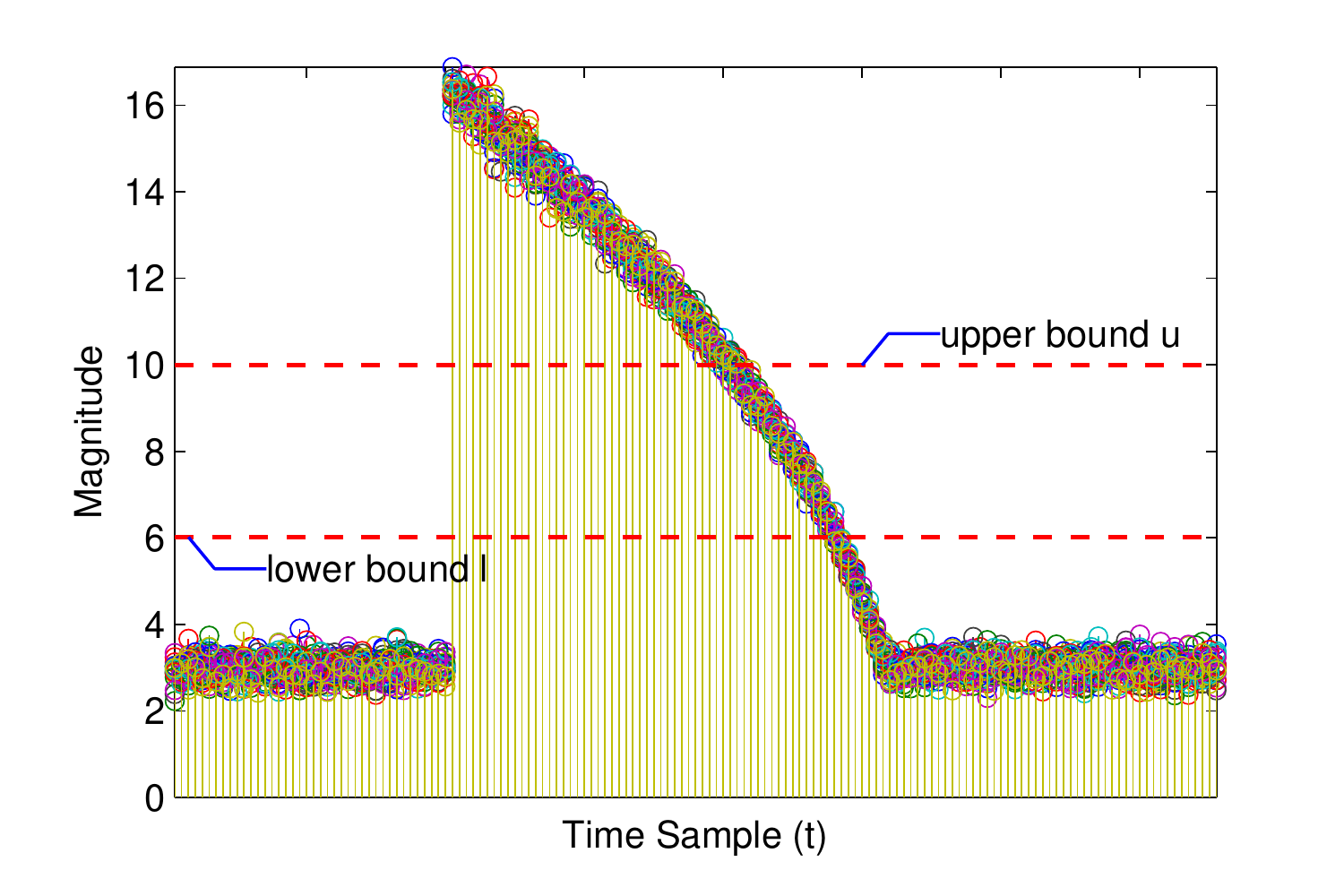}
  \end{center}
\end{minipage}
\caption{Illustration of the performance of the provided bounds of Theorems~\ref{theo:lowerBnd_sigma1_whenNoAttack} and~\ref{theo:upperBnd_sigma1_whenAttack}. A $4$-sparse unobservable attack with $\|\vc{a}\|_2 = 2$ has occurred at $t_a = 129$. Plots depict 300 iterations at each sample time. (a) A window size of $w=8$ and a noise level of $\nu = 0.01$ are considered. (b) A window size of $w=64$ and a noise level of $\nu = 0.04$ are considered.
}
\label{fig:bounds_on_svd_Bus39_UA202_performance}
\end{figure}

Fig.~\ref{fig:bounds_on_svd_Bus39_UA202_performance} illustrates how tight are the probability tail bounds of Theorems~\ref{theo:lowerBnd_sigma1_whenNoAttack} and~\ref{theo:upperBnd_sigma1_whenAttack}.
As mentioned earlier, a $4$-sparse unobservable attack with $\|\vc{a}\|_2 = 2$ has occurred at $t_a = 129$. 
With a given number of measurements $M=85$, we choose $\tau = 4$ and $\epsilon = 0.75$ to make the probability term $2\exp{(-\frac{\tau^2}{2})} + \big((1+\epsilon)e^{-\epsilon}\big)^{M/2} = 6.7\times 10^{-4}$. We stick to these values in all simulations provided in this section. However, different values of $\tau$ and $\epsilon$ can be chosen to achieve a desired tail probability.
A window size of $w=8$ and a noise level of $\nu = 0.01$ are considered in Fig.~\ref{fig:bounds_on_svd_Bus39_UA202_performance}(a). We repeat the simulations for 300 realizations at each sample time.
A window size of $w=64$ and a noise level of $\nu = 0.04$ are considered in Fig.~\ref{fig:bounds_on_svd_Bus39_UA202_performance}(b). As can be seen, the gap between the provided bounds and the actual $\sigma_1(\Delta^t)$ magnitude is larger for cases with larger noise levels. 
%%%%%%%%%%%%%%%%%%%%%%%%%%

\begin{figure}[tb]
\begin{minipage}[b]{\linewidth}
\begin{center}
(a)
 \includegraphics[width =\columnwidth]{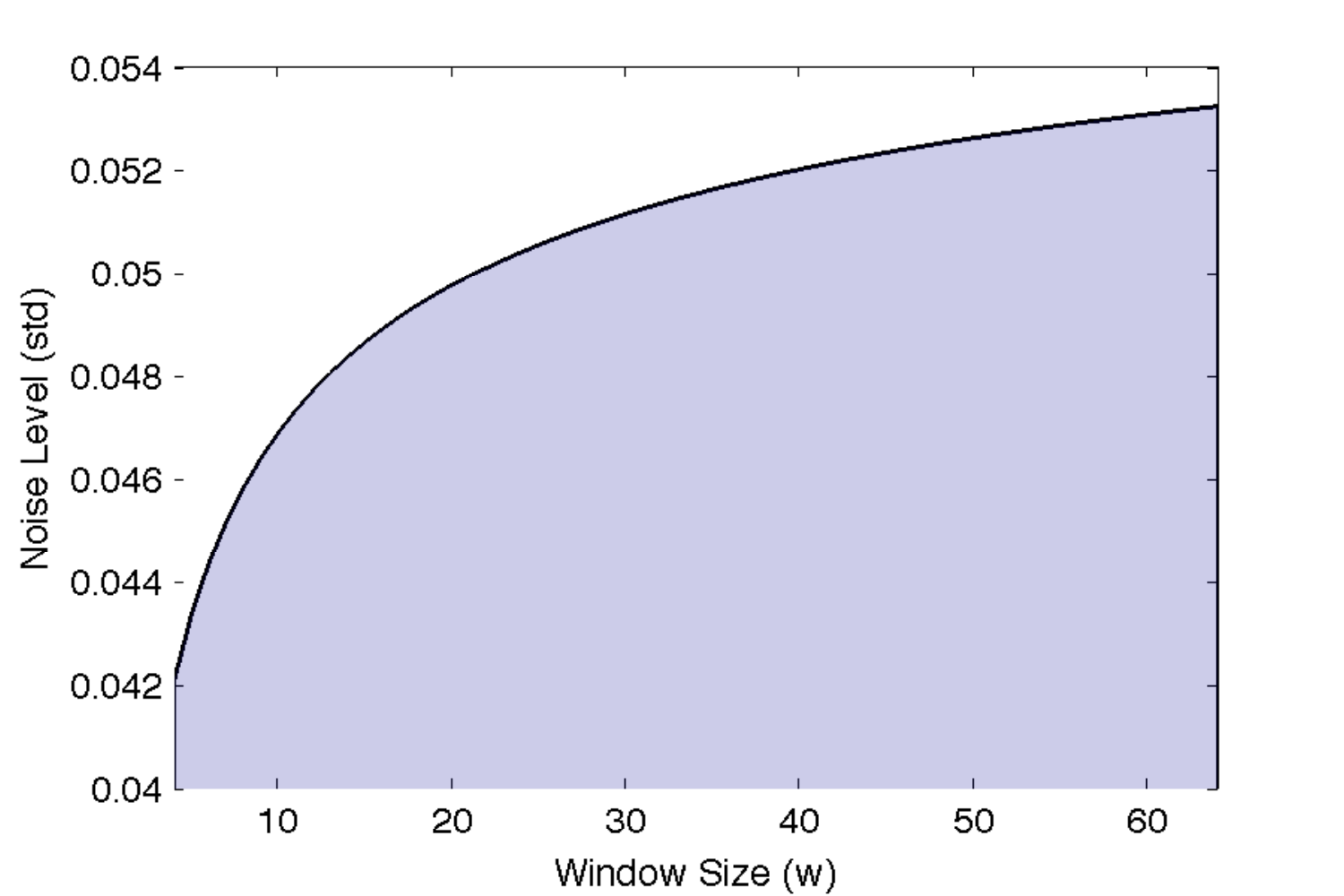}
 \end{center}
\end{minipage}
\begin{minipage}[b]{\linewidth}
\begin{center}
(b)
 \includegraphics[width = \columnwidth]{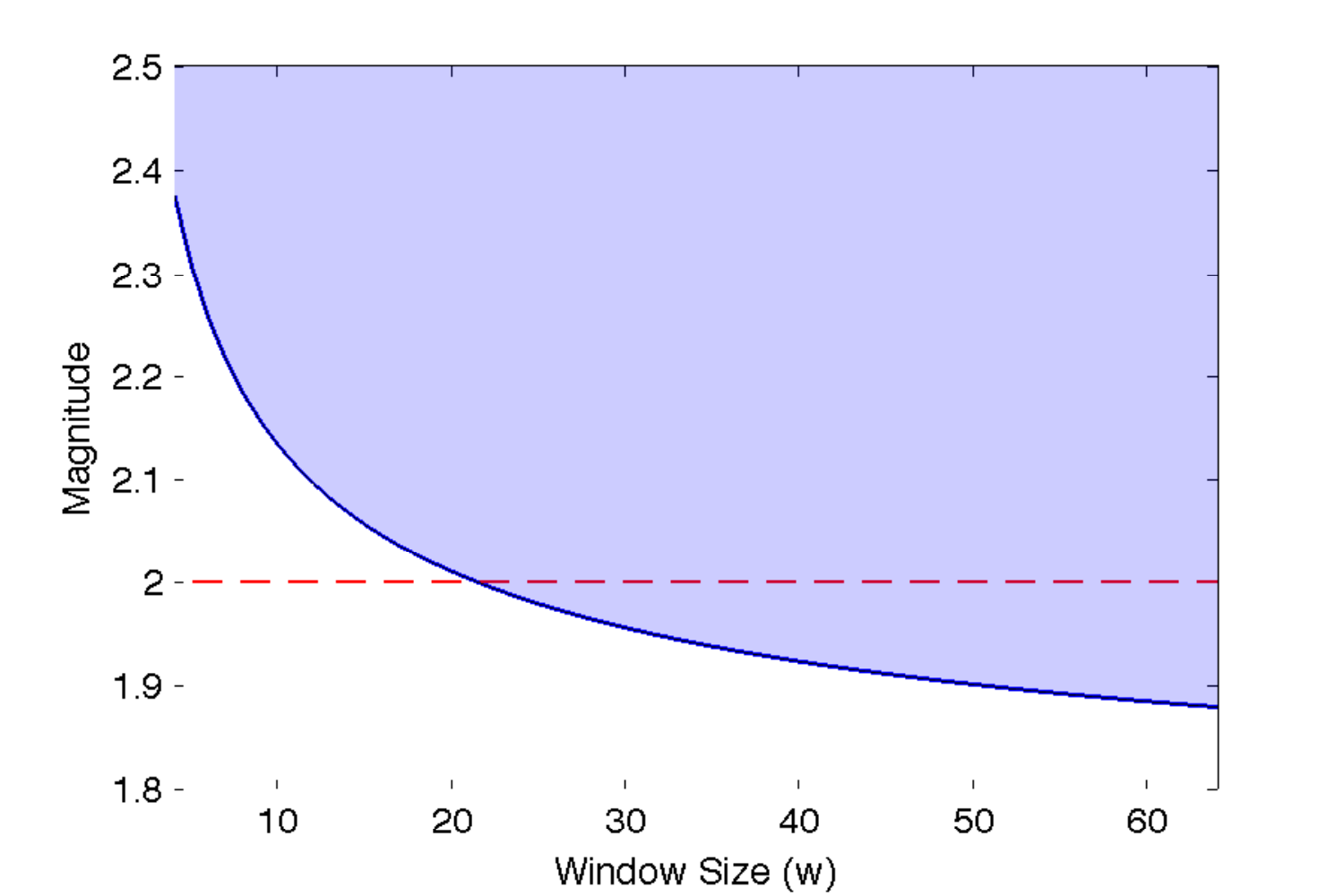}
  \end{center}
\end{minipage}
\caption{Illustration of how different parameters affect the detectability condition of Theorem~\ref{theo:detectionProb}.
(a) $\|\vc{a}\|_2 = 2$ is fixed. One needs to increase $w$ as noise level $\nu$ increases. (b) $\nu = 0.05$ is fixed. For any $\|\vc{a}\|_2$, the curve determines the minimum required window size for having detectability. Attacks with larger magnitudes may be detected with smaller window sizes.}
\label{fig:relation_w_nu_a_Bus39_UA202}
\end{figure}

%%%%%%%%%%%%%%%%%%%%%%%%%%
%%%%%%%%%%%%%%%%%%%%%%%%%%%
Fig.~\ref{fig:relation_w_nu_a_Bus39_UA202} shows how different parameters affect the detectability condition proposed in Theorem~\ref{theo:detectionProb}. We first assume $\|\vc{a}\| = 2$, $M = 85$, $\tau = 4$, and $\epsilon = 0.75$, and we are interested in the relation between window size $w$ and noise level $\nu$ such that the sufficient condition of Theorem~\ref{theo:detectionProb} is satisfied. Fig.~\ref{fig:relation_w_nu_a_Bus39_UA202}(a) shows the result. As can be seen, for larger values of $\nu$, larger window sizes $w$ should be considered such that the attack can be detected with exponentially high probability. In another scenario, we consider the case where $\nu = 0.05$, $M = 85$, $\tau = 4$, and $\epsilon = 0.75$ are fixed and we are interested in finding the relation between window size $w$ and $\|\vc{a}\|_2$. 
Fig.~\ref{fig:relation_w_nu_a_Bus39_UA202}(b) shows the result. For any $\|\vc{a}\|_2$, the curve determines the minimum required window size for having detectability. For example, when $\|\vc{a}\|_2 = 2$ one need to construct the history matrix $\Delta^t$ with $w \geq 22$ such that the sufficient detectability condition of Theorem~\ref{theo:detectionProb} is satisfied. In other words, attacks with larger magnitude may be detected with smaller window sizes and similarly, one one needs to increase the window size as $\|\vc{a}\|_2$ gets smaller.    
%%%%%%%%%%%%%%%%%%%%%%%%%%%

%%%%%%%%%%%%%%%%%%%%%%%%%%
%%%%%%%%%%%%%%%%%%%%%%%%%%%%%%%%%%%%%%%%%%%%%%%%%%%%%%%%%%%%%%%%%%%%%%%%%%%%%%%%

%%%%%%%%%%%%%%%%%%%%%%%%%%%%%%%%%%%%%%%%%%%%%%%%%%%%%%%%%%%%%%%%%%%%%%%%%%%%%%%
%\begin{figure}[tb]
%\begin{minipage}[b]{\linewidth}
%\begin{center}
%(a)
% \includegraphics[width =\columnwidth]{my_fig/detection_Bus39_UA202_w64_q2_q5_test2_plot1}
% \end{center}
% \label{fig:detection_Bus39_UA202_w64_q2_q5_test2_plot1}
%\end{minipage}
%\begin{minipage}[b]{\linewidth}
%\begin{center}
%(b)
% \includegraphics[width = \columnwidth]{my_fig/detection_Bus39_UA202_w64_w8_q2_test2_plot1}
%  \end{center}
%  \label{fig:detection_Bus39_UA202_w64_w8_q2_test2_plot1}
%\end{minipage}
%\caption{Detection performance. For each noise level, the detection rate is calculated over $300$ trials. The number of measurements is $M=85$. (a) The effect of $q$ on the detection performance when $w=64$. (b) The effect of $w$ on the detection performance when $q=2$.}
%\label{fig:detection_performance}
%\end{figure}
%%%%%%%%%%%%%%%%%%%%%%%%%%%%%%%%%%%%%%%%%%%%%%%%%%%%%%%%%%%%%%%%%%%%%%%%%%%%%%%%
%%%%%%%%%%%%%%%%%%%%%%%%%%%%%%%%%%%%%%%%%%%%%%%%%%%%%%%%%%%%%%%%%%%%%%%%%%%%%%%%
%%%%%%%%%%%%%%%%%%%%%%%%%%%%%%%%%%%%%%%%%%%%%%%%%%%%%%%%%%%%%%%%%%%%%%%%%%%%%%%%
    
%%%%%%%%%%%%%%%%%%%%%%%%%%%%%%%%%%%%%%%%%%%%%%%%%%%%%%%%%%%%%%%%%%%%%%%%%%%%%%%%
\appendix
{\textbf{Proof of Theorem~\ref{theo:lowerBnd_sigma1_whenNoAttack}}}
Considering the measurement model given in (\ref{eq:DCmodel-with-attack}), for $t < t_a$ (i.e., before the attack), we have
\begin{equation*}
\Delta^t = G^t + E^t + HX^t, \qquad (t<t_a).
%\label{eq:delta_before_attack}
\end{equation*}
We are interested in showing that there exists an $\ell$ such that $\Prob{\sigma_1(\Delta^t)\geq \ell}$ is very small for $\forall t < t_a$. We have
\begin{align}
\sigma_1(\Delta^t) &= \sigma_1(E^t + G^t + HX^t) \notag \\
&\leq \sigma_1(E^t) + \sigma_1(G^t) + \sigma_1(HX^t) \notag\\
&\leq \sqrt{w}\|\vc{e}^t\|_2+\sigma_1(G^t) + \gamma\sqrt{w}\|H\|,
\label{eq:upperBnd_sigma1D_noAttack}
\end{align}
where we used the assumption that $\|\vc{x}^t-\vc{x}^{t-t_0}\| \leq \gamma$, for all $t \ \text{and} \ t_0 \in \{t_i, \dots,t_f\}$.
Given $\nu, \tau, \epsilon, M, w$, and any $t<t_a$, let's define event $A$ as
\[
\mathcal{E}(A) := \big\{\sigma_1(G^{t}) < \nu(\sqrt{M} +\sqrt{w}+
\tau)\big\}
\]
and event $B$ as
\[
\mathcal{E}(B) := \big \{\|\vc{e}^{t}\|_2 < \nu \sqrt{M}(1+\epsilon)\big\}.
\]
It is trivial to see that if events $\mathcal{E}(A)$ and $\mathcal{E}(B)$ happen, then event $C$ defined as  
\begin{equation*}
\mathcal{E}(C) := 
\big \{  \sqrt{w}\|\vc{e}^t\|_2 +\sigma_1(G^t) + \gamma\sqrt{w}\|H\|< \ell \big \}
\end{equation*}
happens where
\[
\ell := \nu\sqrt{w}\sqrt{M}(1+\epsilon)+\nu(\sqrt{M}+\sqrt{w}+\tau)+\gamma\sqrt{w}\|H\|.
\]
Using results from \ac{CoM} phenomenon of random processes~\cite{ledoux2001concentration,lugosi2004concentration,sanandaji2013com}, we first show that for any $t$, random variables $\sigma_1(G^{t})$ and $\|\vc{e}^{t}\|_2$ are highly concentrated around their expected value. The following lemmas provide such \ac{CoM} bounds. 
\begin{lemma}(\cite{sanandaji2013com},\cite[Lemma 2]{sanandaji2012tutorial})
Let $\vc{e}$ be a vector in $\real^{M}$ whose entries are independent Gaussian random variables with zero mean and $\nu^2$ variance. Then for every $\epsilon \geq 0$, 
\[
\Prob{\|\vc{e}\|_2 \geq \nu \sqrt{M}(1+\epsilon)} \leq \big((1+\epsilon)e^{-\epsilon}\big)^{M/2}.
\]
\label{lem:norm_1_conc}
\end{lemma}
\begin{lemma}(\cite{vershynin2011introduction})
Let $G$ be an $M \times w$ matrix whose entries are independent Gaussian random variables with zero mean and $\nu^2$ variance. Then for every $\tau \geq 0$, 
\[
\Prob{\sigma_1(G) \geq \nu(\sqrt{M} +\sqrt{w}+\tau)}\leq 2\exp{(-\frac{\tau^2}{2})}.
\]
\label{lem:delta_1_conc}
\end{lemma}
Using Lemma~\ref{lem:norm_1_conc} and Lemma~\ref{lem:delta_1_conc}, and noting that $\Prob{\mathcal{E}(C)^c}\leq \Prob{\mathcal{E}(A)^c}+\Prob{\mathcal{E}(B)^c}$, we have
\begin{align*}
\Prob{\sigma_1(\Delta^t)\geq \ell} &\leq \Prob{\sqrt{w}\|\vc{e}^t\|_2+\sigma_1(G^t) \geq \ell} \\ &\leq 2\exp{(-\frac{\tau^2}{2})} + \big((1+\epsilon)e^{-\epsilon}\big)^{M/2},
\end{align*}
where we used (\ref{eq:upperBnd_sigma1D_noAttack}) in showing the first inequality.
\hfill $\blacksquare$
%%%%%%%%%%%%%%%%%%%%%%%%%%%%%%%%%%%%%%%%%%%%%%%%%%%%%%%%%%%%%%%%%%%%%%%%%%%%%%%%
%{\textbf{Proof of Lemma~\ref{lem:bounds_on_svd_1}}}
%
%Modifying~\cite[Theorem 6]{merikoski2004inequalities}, we have
%%
%\[
%\sigma_i(\Delta^{t_a}) = \sigma_i(G^{t_a}+E^{t_a}+A) \leq \sigma_1(G^{t_a})+\sigma_i(E^{t_a}+A),
%\]
%for any $i \geq 1$. In particular, for $i=1$
%%
%\begin{equation*}
%\sigma_1(\Delta^{t_a}) \leq \sigma_1(G^{t_a})+\sigma_1(E^{t_a}+A).
%\end{equation*}
%%
%Noting that $E^{t_a}+A$ is an $M \times w$ rank-$1$ matrix, it is trivial to see that $\sigma_{1}(E^{t_a}+A) = \sqrt{w}\|\vc{e}^{t_a}+\vc{a}\|_2$. 
%%
%Therefore,
%\begin{equation}
%\sigma_1(\Delta^{t_a}) \leq \sigma_1(G^{t_a}) + \sqrt{w}\|\vc{e}^{t_a}+\vc{a}\|_2.
%%\label{eq:delta_2_upper_bnd}
%\end{equation}
%
%Similarly, we have
%%
%\[
%\sigma_i(\Delta^{t_a-1}) = \sigma_i(G^{t_a-1}+E^{t_a-1}) \leq \sigma_1(G^{t_a-1})+\sigma_i(E^{t_a-1}),
%\]
%for any $i \geq 1$. In particular, for $i=1$
%%
%\begin{align}
%\sigma_1(\Delta^{t_a-1}) &\leq \sigma_1(G^{t_a-1})+\sigma_1(E^{t_a-1}) \notag\\
%&= \sigma_1(G^{t_a-1}) + \sqrt{w}\|\vc{e}^{t_a-1}\|_2.
%\label{eq:delta_1_upper_bnd}
%\end{align}
%

%%%%%%%%%%%%%%%%%%%%%%%%%%%%%%%%%%%%%%%%%%%%%%%%%%%%%%%%%%%%%%%%%%%%%%%%%%%%%
\textbf{Proof of Theorem~\ref{theo:upperBnd_sigma1_whenAttack}}
First we need a lower bound on $\sigma_1(\Delta^{t_a})$. Modifying~\cite[Theorem 6]{merikoski2004inequalities}, we can derive a lower bound on the first singular value of $\Delta^{t_a}$ as
\begin{align}
\sigma_1(\Delta^{t_a}) &= \sigma_1(E^{t_a}+G^{t_a}+HX^{t_a}+A) \notag \\
&\geq \big\vert \sigma_{i}(E^{t_a}+A) - \sigma_{i}(G^{t_a}+HX^{t_a})\big\vert
\end{align}
for all $1 \leq i \leq \min\{M,w\}$.
In particular, for $i=1$
\begin{equation}
\sigma_1(\Delta^{t_a}) \geq \big\vert \sigma_{1}(E^{t_a}+A) - \sigma_{1}(G^{t_a}+HX^{t_a}) \big\vert.
\label{eq:lower_bnd_delta_1}
\end{equation}
%
%\begin{lemma}
%Let $A$ be an $N \times w$ matrix containing $w$ equal columns $\vc{a} \in \real^N$. The maximum singular value of $A$ is equal to $\sqrt{w}\|\vc{a}\|_2$.
%\label{lem:svd1}
%\end{lemma}
%
%%
%\textbf{Proof of Lemma~\ref{lem:svd1}}
%For any matrix $A \in \real^{N \times w}$ of rank $r$, we have
%\[
%\|A\|_2 \leq \|A\|_F \leq \sqrt{r}\|A\|_2.
%\]
%%
%Noting that a matrix of equal columns is a rank-1 matrix, the result follows.
%\hfill $\blacksquare$

%
Assuming $\sigma_{1}(E^{t_a}+A) > \sigma_{1}(G^{t_a}+HX^{t_a})$, we have
%then 
%
\begin{align*}
\sigma_1(\Delta^{t_a}) &\geq 
\sqrt{w}\|\vc{e}^{t_a}+\vc{a}\|_2 - \sigma_{1}(G^{t_a}+HX^{t_a})\\
& \geq \sqrt{w}\|\vc{e}^{t_a}+\vc{a}\|_2 - \sigma_{1}(G^{t_a}) - \gamma\sqrt{w}\|H\|,
%\label{eq:delta_1_lower_bnd}
\end{align*}
%\hfill $\blacksquare$
where we used the fact that $\sigma_{1}(E^{t_a}+A)  = \sqrt{w}\|\vc{e}^{t_a}+\vc{a}\|_2$ and $\sigma_{1}(G^{t_a}+HX^{t_a}) \leq \sigma_{1}(G^{t_a})+\sigma_1(HX^{t_a})$.
Assuming $\|\vc{a}\|_2 \geq \|\vc{e}^{t_a}\|_2$ and using the reverse triangle inequality, 
\begin{equation}
\|\vc{e}^{t_a}+\vc{a}\|_2 \geq \big | \|\vc{a}\|_2 - \|\vc{e}^{t_a}\|_2 \big| = \|\vc{a}\|_2 - \|\vc{e}^{t_a}\|_2.  
\label{eq:reverse_tri_ineq_1}
\end{equation}
Therefore, from~(\ref{eq:lower_bnd_delta_1})
\begin{align}
\label{eq:lowerBnd_sigma1_withAttack}
\sigma_1(\Delta^{t_a}) &\geq \sqrt{w}\|\vc{e}^{t_a}+\vc{a}\|_2-\sigma_1(G^{t_a})- \gamma\sqrt{w}\|H\| \\
&\geq \sqrt{w}\|\vc{a}\|_2 - \sqrt{w}\|\vc{e}^{t_a}\|_2-\sigma_1(G^{t_a})- \gamma\sqrt{w}\|H\|.\notag
\end{align}

Given $\nu, \tau, \epsilon, M, w$, and $t_a$, let's define event $A$ as
\[
\mathcal{E}(A) := \big\{\sigma_1(G^{t_a}) < \nu(\sqrt{M} +\sqrt{w}+
\tau)\big\}
\]
and event $B$ as $\mathcal{E}(B) := \big \{\|\vc{e}^{t_a}\|_2 < \nu \sqrt{M}(1+\epsilon)\big\}$.
It is trivial to see that if events $\mathcal{E}(A)$ and $\mathcal{E}(B)$ happen, then event $C$ defined as  
\begin{equation*}
\mathcal{E}(C) := \\
\big \{\sqrt{w}\|\vc{a}\|_2 - \sqrt{w}\|\vc{e}^{t_a}\|_2-\sigma_1(G^{t_a}) - \gamma\sqrt{w}\|H\|> u\big \}
\end{equation*}
happens where $u := \sqrt{w}\|\vc{a}\|_2 - \ell$ and
\[
\ell = \nu\sqrt{w}\sqrt{M}(1+\epsilon)+\nu(\sqrt{M}+\sqrt{w}+\tau)+\gamma\sqrt{w}\|H\|.
\]
Using Lemma~\ref{lem:norm_1_conc} and Lemma~\ref{lem:delta_1_conc}, and noting that $\Prob{\mathcal{E}(C)^c}\leq \Prob{\mathcal{E}(A)^c}+\Prob{\mathcal{E}(B)^c}$, we have
\begin{multline*}
\Prob{\sigma_1(\Delta^{t_a}) \leq  u} \leq \\ \Prob{\sqrt{w}\|\vc{a}\|_2 - \sqrt{w}\|\vc{e}^{t_a}\|_2-\sigma_1(G^{t_a}) \leq u} \\
\leq 2\exp{(-\frac{\tau^2}{2})} + \big((1+\epsilon)e^{-\epsilon}\big)^{M/2},
\label{eq:upperBnd_sigma1_whenAttack}
\end{multline*}
where 
we used (\ref{eq:lowerBnd_sigma1_withAttack}) in showing the first inequality. 
\hfill $\blacksquare$

%%%%%%%%%%%%%%%%%%%%%%%%%%%%%%%%%%%%%%%%%%%%%%%%%%%%%%%%%%%%%%%%%%%%%%%%%%%%%
\textbf{Proof of Theorem~\ref{theo:detectionProb}}
For any $\ell$ and $u$ such that $\ell<u$, 
\[
\Prob{\text{detection}} = 
\Prob{\sigma_1(\Delta^t)<\ell \ \text{and} \
\sigma_1(\Delta^{t_a})>u}.
\]
Using similar techniques as used in the proof of 
Theorems~\ref{theo:lowerBnd_sigma1_whenNoAttack} and \ref{theo:upperBnd_sigma1_whenAttack}, we have
\begin{align*}
\Prob{\text{not detection}} &\leq 
\Prob{\sigma_1(\Delta^t)\geq \ell} +  \Prob{\sigma_1(\Delta^{t_a})\leq u} \\
&\leq 2\big[2\exp{(-\frac{\tau^2}{2})} + \big((1+\epsilon)e^{-\epsilon}\big)^{M/2}\big],
\end{align*}
where $\ell$ and $u$ are chosen as given in Theorems~\ref{theo:lowerBnd_sigma1_whenNoAttack} and \ref{theo:upperBnd_sigma1_whenAttack}, respectively. If 
\begin{equation*}
\|\vc{a}\|_2 > 2\big[\nu\sqrt{M}(1+\epsilon+
\frac{1}{\sqrt{w}}+\frac{1}{\sqrt{M}}+\frac{\tau}{\sqrt{M}\sqrt{w}})+\gamma\|H\|_2\big],
\end{equation*}
then $\ell<u$ and this completes the proof.
\hfill $\blacksquare$

%
%\input appendix.tex

%{\textbf{Proof of Lemma~\ref{lem:trivial_1}}}
%The proof is trivial and is a based on the triangle inequality. We have
%\begin{align*}
%\|\vc{a}\|_2 &\geq \|K\|_2\|H\|_2Q\\
%\|H\|_2\|\vc{c}\|_2 \geq \|H\vc{c}\|_2 &\geq \|K\|_2\|H\|_2Q\\
%\|\vc{c}\|_2 &\geq \|K\|_2Q,
%\end{align*}
%where $\|H\|_2$ and $\|K\|_2$ are spectral norms of $H$ and $K$, and $\|\vc{a}\|_2$ and $\|\vc{c}\|_2$ are Euclidean norms of $\vc{a}$ and $\vc{c}$.\hfill $\blacksquare$

\bibliographystyle{IEEEtran}
\bibliography{AttackDetection.bib}

\end{document}